\documentclass[11pt]{article}
\usepackage[latin1]{inputenc}
\usepackage{amsmath,amsfonts,amssymb}
\usepackage{amscd}
\usepackage{graphicx}

\usepackage[T1]{fontenc}

\newcommand{\al}{\alpha}
\newcommand{\be}{\beta}
\newcommand{\ga}{\gamma}

\newcommand{\proof}{\noindent{\bf Proof. }}
\newcommand{\qed}{{$\diamond$ \\}}

\newcommand{\Z}{{\mathbb{Z}}}
\newcommand{\N}{{\mathbb{N}}}

\newcommand{\End}{\operatorname{End}}
\newcommand{\Hom}{\operatorname{Hom}}

\newcommand{\add}{\operatorname{add}}
\newcommand{\HH}{\operatorname{HH}}

\renewcommand{\P}{{\cal P}}

\newtheorem{theorem}{Theorem}[section]

\newtheorem{cor}[theorem]{Corollary}
\newtheorem{lemma}[theorem]{Lemma}
\newtheorem{example}[theorem]{Example}
\newtheorem{remark}[theorem]{Remark}

\begin{document}
\begin{center}

{\LARGE\bf
$q$-Cartan matrices and combinatorial
invariants of derived categories
for skewed-gentle algebras}

\bigskip

\bigskip

{\Large\bf Christine Bessenrodt}
\footnote{e-mail: {\tt bessen@math.uni-hannover.de}}
\smallskip

{\large Fakult\"at f\"ur Mathematik und Physik, Universit\"at Hannover,\\
Welfengarten 1, D-30167 Hannover, Germany}

\bigskip

{\Large\bf Thorsten Holm} \footnote{e-mail: {\tt tholm@maths.leeds.ac.uk}}
\smallskip

{\large Department of Pure Mathematics, University of Leeds,\\
Leeds LS2 9JT, England}

%\medskip

%June 25, 2005-version, revised version: Sept 22, 2005

\bigskip

%MSC: 16G20,18E30,15A15, 05E99,05C50

\begin{abstract}
Cartan matrices are of fundamental importance in representation theory.
For algebras defined by quivers
with monomial relations
the computation of the entries of the Cartan matrix amounts to
counting nonzero paths in the quivers, leading naturally
to a combinatorial setting.
Our main motivation are derived module categories and their invariants:
the invariant factors, and hence the determinant, of the Cartan
matrix are preserved by derived equivalences.

The paper deals with the class of (skewed-)
gentle algebras which occur naturally in representation theory,
especially in the context of derived categories.
%These algebras are defined in purely combinatorial terms.
We study $q$-Cartan matrices, where each nonzero path is
weighted by a power of an indeterminate~$q$ according to its length.
Specializing $q=1$ gives the classical Cartan matrix.
We determine normal forms for the $q$-Cartan matrices of
skewed-gentle algebras. In particular, we give explicit
combinatorial formulae for the
invariant factors and thus also for the determinant.
%This generalizes the results of~\cite{TH}.
As an application of our main results we show how to use
our formulae for the difficult problem of
distinguishing derived equivalence classes.
\medskip

\noindent
MSC-Classification: 16G10, 18E30, 05E99, 05C38, 05C50
%CB: 15A36? 16G60?
%vielleicht eher noch 05E99, 05C38 (Paths and cycles, habe ich gerade erst entdeckt) oder 05C50
\end{abstract}

\end{center}

\section{Introduction} \label{Sec-intro}

This paper deals with combinatorial aspects in the
representation theory of algebras. More precisely,
for certain classes of algebras which are defined purely combinatorially
by directed graphs
and homogeneous relations
we will characterize important
represen\-ta\-tion-theoretic invariants in a combinatorial
way. In particular, this leads to new explicit invariants for
the derived module categories of the algebras involved.
%Derived categories are nowadays one of the most important
%objects in representation theory and in geometry, with numerous
%applications, e.g.\ in mathematical physics.
%In the representation theory of algebras, a difficult
%problem is to find invariants of the algebras which are preserved under
%derived equivalences.
%Although quite a few general invariants of derived module
%categories are known,
%most of them are very hard to compute.
%This makes it usually difficult to distinguish algebras up
%to derived equivalence.

%The invariants we consider in this paper are the invariant factors of the
%Cartan matrix. Cartan matrices play a fundamental r\^ole in all areas
%of representation theory, containing complete information on the
%multiplicities of simple modules as composition factors of
%projective indecomposable modules. The entries of the Cartan matrix
%are not known in general and some famous open problems are related
%to their computation.

The starting point for this article is that the unimodular equivalence
class of the Cartan matrix of a finite dimensional algebra is
invariant under derived equivalence.
Hence, being able to determine normal forms of Cartan matrices
yields invariants of the derived category.

The class of algebras we study are the gentle algebras, and the related
skewed-gentle algebras.
Gentle algebras are defined purely
combinatorially in terms of a quiver with relations
(for details, see Section~\ref{Sec-def});
the more general skewed-gentle algebras
(introduced in~\cite{GdlP})
are then defined from gentle algebras
by specifying special vertices which are split for the
quiver of the skewed-gentle algebra (see Section~\ref{Sec-skewed}).
These algebras occur naturally in the
representation theory of finite dimensional algebras, especially in the context of
derived categories. For instance, the algebras which are derived
equivalent to hereditary algebras of type~$\mathbb{A}$ are precisely
the gentle algebras whose underlying undirected graph is a tree~\cite{AH}.
The algebras which are derived equivalent to hereditary
algebras of type~$\tilde{\mathbb{A}}$ are certain gentle
algebras whose underlying graph has exactly one cycle~\cite{AS}.
Remarkably, the class of gentle algebras is closed under derived
equivalence~\cite{SZ};
but note that the class of skewed-gentle algebras is not
closed under derived equivalence.

A fundamental distinction in the representation theory
of algebras is given by the representation type, which can be
either finite, tame or wild.
In the modern context of derived categories, also derived
representation types have been defined.
Again, gentle algebras occur naturally in this context.
D.~Vossieck~\cite{V} showed that an algebra~$A$
has a discrete derived category
if and only if either $A$ is derived equivalent to a hereditary
algebra of type $\mathbb{A}$, $\mathbb{D}$, $\mathbb{E}$ or
$A$ is gentle with underlying quiver $(Q,I)$ having exactly one
(undirected) cycle
and the number of clockwise and of counterclockwise paths of
length~2 in the cycle that belong to~$I$ are different.
Skewed-gentle algebras are known to be of derived tame representation
type (for a definition of derived tameness, see~\cite{GK}).

It is a long-standing open problem to classify gentle algebras
up to derived equivalence. A complete answer has only been obtained
for the derived discrete case~\cite{BGS}. The main problem is to find
good invariants of the derived categories.

In this paper we provide easy-to-compute invariants of the
derived categories of skewed-gentle algebras which are of
a purely combinatorial nature. Our results are obtained from a
detailed computation of the $q$-Cartan matrices of gentle
and skewed-gentle algebras, respectively.

The following notion will be crucial throughout the paper. Let
$(Q,I)$ be a (gentle) quiver with relations. An oriented path
$p=p_0p_1\ldots p_{k-1}$ with arrows $p_0,\ldots,p_{k-1}$ in~$Q$
is called an oriented $k$-cycle
with full zero relations if $p$ has the same start and end point,
and if $p_ip_{i+1}\in I$ for all $i=0,\ldots,k-2$
and also $p_{k-1}p_0\in I$. Such a cycle is called
minimal if the arrows $p_0,p_1,\ldots,p_{k-1}$
on~$p$ are pairwise different.

We call two matrices $C,D$ with entries in a polynomial
ring~$\mathbb{Z}[q]$
unimodularly equivalent (over~$\mathbb{Z}[q]$)
if there exist matrices~$P,Q$
over~$\mathbb{Z}[q]$ of determinant~1 such that $D=PCQ$.

We can now state our main result on gentle algebras.
\smallskip

{\bf Theorem~1.} \label{Thm1} {\em
Let $(Q,I)$ be a gentle quiver, and $A= KQ/I$ the
corresponding gentle algebra.
Denote by $c_k$ the number of minimal oriented
$k$-cycles in~$Q$ with full zero relations.\\
Then the $q$-Cartan matrix~$C_A(q)$
is unimodularly equivalent (over~$\Z[q]$) to a diagonal matrix with entries
$(1-(-q)^k)$, with multiplicity ${c_k}$, $k\geq 1$,
and all further diagonal entries being~$1$.
}
\medskip

This theorem has the following direct consequences.
\smallskip

{\bf Corollary 1.} {\em
Let $(Q,I)$ be a gentle quiver, and
$A=KQ/I$ the corresponding gentle algebra.
Denote by $c_k$ the number of minimal
oriented $k$-cycles in~$Q$
with full zero relations.
Then the $q$-Cartan matrix~$C_A(q)$
has determinant $$\det C_A(q)= \prod_{k \geq 1} (1-(-q)^k)^{c_k} \:.$$
}
\medskip

The following consequence of Corollary~1 was first proved in~\cite{TH}.
For a gentle quiver~$(Q,I)$
we denote
by~$oc(Q,I)$ the number of minimal oriented cycles of odd length
in~$Q$ having full zero relations, and by~$ec(Q,I)$ the number
of analogous cycles of even length.
\smallskip

{\bf Corollary 2.}
{\em
Let $(Q,I)$ be a gentle quiver, and $A=KQ/I$ the corresponding
gentle algebra.
Then for the determinant of the Cartan
matrix~$C_A$ the following holds.
$$\det C_A = \left\{
\begin{array}{ll} 0 & \textrm{if }\: ec(Q,I)>0 \\
                  2^{oc(Q,I)} & \textrm{else}
\end{array}
\right.
$$
}

The most important application of Theorem~1
is the following corollary which gives
for gentle algebras  easy-to-check  combinatorial invariants
of the derived category.
\medskip

{\bf Corollary 3.} {\em
Let $(Q,I)$ and $(Q',I')$ be gentle quivers,
and let $A=KQ/I$ and $A'=KQ'/I'$ be the corresponding
gentle algebras.
If $A$ and~$A'$ are derived equivalent, then
$ec(Q,I)=ec(Q',I')$ and $oc(Q,I)=oc(Q',I')$.
}
\medskip

As an illustration we give in Section~\ref{sec:Gentles} a
complete derived equivalence
classification of gentle algebras with two simple modules
and of gentle algebras with three simple modules and Cartan
determinant~0.
\medskip

Our main result on skewed-gentle algebras determines the normal form of
their $q$-Cartan matrices.
\smallskip

{\bf Theorem 2.} {\em
Let $\hat A=K\hat Q/\hat I$ be a skewed-gentle algebra, arising from choosing
a suitable set of special vertices
in the gentle quiver~$(Q,I)$.
Denote by $c_k$ the number of minimal oriented $k$-cycles in~$Q$ with full
zero relations.\\
Then the $q$-Cartan matrix~$C_{\hat A}(q)$ is unimodularly equivalent
to a diagonal matrix with entries
$1-(-q)^k$, with multiplicity ${c_k}$, $k\geq 1$,
and all further diagonal entries being~$1$.
}
\medskip

As an immediate consequence we obtain that the Cartan determinant
of any skewed-gentle algebra is the same as the Cartan determinant
for the underlying gentle algebra.
\medskip

{\bf Corollary 4.} {\em
Let $\hat A=K\hat Q/\hat I$ be a skewed-gentle algebra, arising
from choosing a suitable set of special vertices in the gentle
quiver~$(Q,I)$, with corresponding gentle
algebra~$A=KQ/I$. Then
$\det C_{\hat A}(q)= \det C_A(q)$, and thus in particular,
the determinants of the ordinary Cartan matrices coincide, i.e.,
$\det C_{\hat A}= \det C_{A}$.
}
\medskip

The paper is organized as follows. In Section~\ref{Sec-def}
we collect the necessary background and definitions about
quivers with relations and ($q$-)Cartan matrices.
In Section~\ref{sec:Gentles} we prove all the main results about
$q$-Cartan matrices for gentle algebras. Here we also give
some extensive examples to illustrate our results.
Section~\ref{Sec-skewed}
contains the analogous main results for skewed-gentle algebras.
\medskip

\noindent
{\bf Acknowledgement.} We thank the referee for
very helpful and insightful comments. In particular, we
are grateful for pointing out the importance of graded
derived equivalences in our context of $q$-Cartan matrices.

%%%%%%%%%%%%%%%%%%%%%%%%%%%%%%%%%%%%%%%%%%%%%%%%%%%%%%%%%%%%%%%%%%%

\section{Quivers, $q$-Cartan matrices and derived invariants}
\label{Sec-def}

Algebras can be defined naturally from a combinatorial setting
by using directed graphs. A finite directed graph~$Q$ is called a
{\em quiver}. For any arrow $\alpha$ in~$Q$ we denote by $s(\alpha)$
its start vertex and by $t(\alpha)$ its end vertex. An oriented
path~$p$ in~$Q$ of length~$r$
is a sequence $p=\alpha_1\alpha_2\ldots\alpha_r$ of arrows
$\alpha_i$ such that $t(\alpha_i)=s(\alpha_{i+1})$ for all
$i=1,\ldots,r-1$.
(Note that for each vertex~$v$ in~$ Q$ we allow
a trivial path~$e_v$ of length~$0$, having~$v$ as its start
and end vertex.)
For such a path~$p$ we then denote by
$s(p):=s(\alpha_1)$ its start vertex and by $t(p):=t(\alpha_r)$
its end vertex.

The path algebra~$KQ$, where $K$ is any field, has as
basis the set of all oriented paths in~$Q$. The multiplication
in the algebra~$KQ$ is defined
by concatenation of paths, i.e.,
the product of two paths $p$ and~$q$ is defined
to be the concatenated path~$pq$ if $t(p)=s(q)$, and zero
otherwise.
%Note that our convention is to write paths from
%left to right.

%Such a path algebra $KQ$ is finite dimensional precisely when
%$Q$ does not contain an oriented cycle.
%Path algebras form an
%important, but rather restricted class of algebras. Far
More general algebras can be obtained by introducing relations on a
path algebra.
An ideal $I\subset KQ$ is called admissible if $I\subseteq
J^2$ where $J$ is the ideal of~$KQ$
generated by the arrows of~$Q$.

The pair $(Q,I)$ where $Q$ is a quiver and $I\subset KQ$
is an admissible ideal is called a {\em quiver with relations}.

%A famous theorem of P.\,Gabriel
%states that if $K$ is algebraically closed, any finite-dimensional
%$K$-algebra is Morita equivalent to a factor algebra
%$KQ/I$ where $I$ is an admissible ideal.

%So for many representation-theoretic purposes
%it suffices to consider algebras of the form
%$KQ/I$, often referred to as {\em quivers with relations}.

For any quiver with relations $(Q,I)$, we can consider the
factor algebra~$A=KQ/I$, where $K$ is any field.
%By a slight abuse of notation we
We identify paths in the quiver~$Q$ with their
cosets in~$A$. Let $Q_0$ denote the set of vertices of~$Q$.
For any $i\in Q_0$ there is a path~$e_i$ of length~zero.
These are primitive orthogonal idempotents in~$A$,
the sum $\sum_{i\in Q_0} e_i$ is the unit element in~$A$.
In particular we get $A=1\cdot A=\oplus_{i\in Q_0} e_iA$,
hence the (right) $A$-modules
$P_i:=e_iA$ are the indecomposable projective
$A$-modules.

The {\em Cartan matrix} $C=(c_{ij})$ of an algebra~$A=KQ/I$
is the $|Q_0|\times |Q_0|$-matrix defined by setting
$c_{ij}:=\dim_K \Hom_A(P_j,P_i)$.

Recall that when  $I$ is
generated by monomials, $A=KQ/I$ is called a monomial algebra.
For monomial algebras,
computing entries of the Cartan matrix reduces to counting
paths in the quiver~$Q$ which are nonzero in~$A$.
In fact,
any homomorphism $\varphi:e_jA\to e_iA$ of right $A$-modules
is uniquely determined by $\varphi(e_j)\in e_iAe_j$,
the $K$-vector space generated by all paths in~$Q$ from vertex~$i$
to vertex~$j$, which are nonzero in~$A=KQ/I$.
In particular, we have $c_{ij}=\dim_Ke_iAe_j$.

This is the key viewpoint in this paper, enabling us to
obtain results on the representation-theoretic Cartan invariants
by combinatorial methods.
%\smallskip
%In this paper we shall
It allows to study a
%the following
refined version of the
Cartan matrix, which we call
the $q$-Cartan matrix. (It also occurred in the literature as
filtered Cartan matrix, see for instance~\cite{F}.)

Let $Q$ be a quiver and assume that
the relation ideal~$I$ is generated by homogeneous relations, i.e., by linear
combinations of paths having the same length (actually, for the
algebras considered in this paper, the ideal~$I$ will always be generated
by monomials and commutativity (mesh) relations).
The path algebra~$KQ$ is a graded algebra, with grading given by
path lengths. Since $I$ is homogeneous, the factor algebra~$A=KQ/I$
inherits this grading. So
the morphism spaces $\Hom_A(P_j,P_i)\cong e_iAe_j$ become
graded vector spaces. Recall that the dimensions of these
vector spaces are the entries of the (ordinary) Cartan matrix.
\medskip

\noindent
{\bf Definition.} Let $A=KQ/I$ be a finite-dimensional algebra,
and assume that the ideal~$I$ is generated by homogeneous relations.
For any vertices~$i$ and~$j$ in~$Q$
let $e_iAe_j= \oplus_{n} (e_iAe_j)_n$ be the
graded components.

Let $q$ be an indeterminate.
The $q$-Cartan matrix $C_A(q)=(c_{ij}(q))$ of~$A$ is defined as the matrix
with entries
$c_{ij}(q):= \sum_n \dim_K (e_iAe_j)_n\, q^n\in \mathbb{Z}[q]$.
\medskip

In other words, the entries of the $q$-Cartan matrix are the
Poincar\'{e} polynomials of the graded homomorphism spaces
between projective modules. Loosely speaking, when counting paths
in the quiver of the algebra, each path is weighted by
some power of~$q$ according to its length.
%CB: expliziter Hinweis auf forthcoming paper, bei dem dies Thema ist?

Clearly, specializing $q=1$ gives back the usual Cartan matrix~$C_A$
(i.e., we forget the grading).
Even if we are mainly interested in the ordinary Cartan matrix, %at $q=1$,
the point of view of~$q$-Cartan matrices
provides some new insights as we take a closer look
at the invariants of the Cartan matrix.

\begin{example} \label{brauer-ex}
{\em We consider the following two quivers.

%%%%%%%%%%%%%%%%%%%%%%

\begin{center}
\includegraphics[scale=0.5]{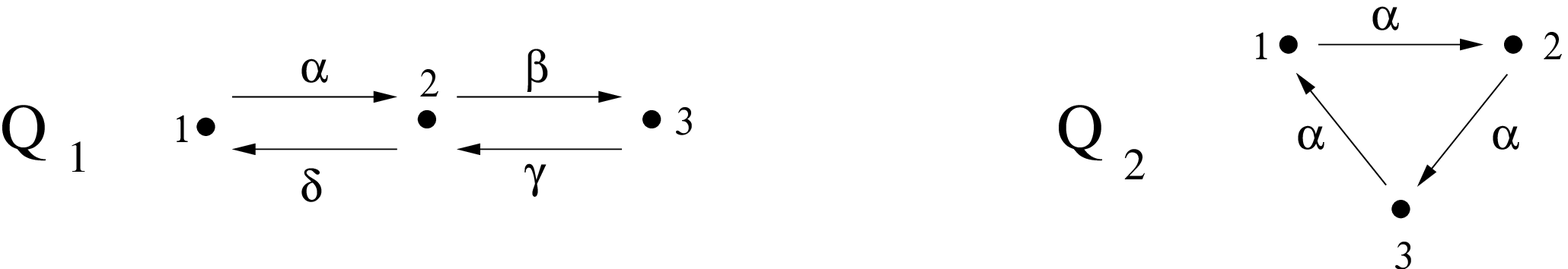}
\end{center}

%%%%%%%%%%%%%%%%%%%%%%

%From these we define the following two algebras $A$ and $B$
%by quivers with relations.
Let $A=KQ_1/I_1$,
%be defined by the quiver $Q_1$,
where the ideal~$I_1$ is generated by $\alpha\beta$,
$\gamma\delta$ and~$\delta\alpha-\beta\gamma$.
The $q$-Cartan matrix of~$A$ has the form
%can be obtained by counting nonzero
%paths in the quiver and weighting them according to their lengths.
%We get
$$C_A(q) = \mbox{{\small $\left( \begin{array}{ccc}
1+q^2 & q & 0 \\ q & 1+q^2 & q \\ 0 & q & 1+q^2
\end{array} \right)$}}.
$$

The second algebra $B=KQ_2/I_2$ is defined by the quiver~$nQ_2$,
subject to the generating relations~$\alpha^4$ (i.e.\ all paths of
length four are zero). The $q$-Cartan matrix of~$B$ has the form
$$C_B(q) = \mbox{{\small $\left( \begin{array}{ccc}
1+q^3 & q & q^2 \\ q^2 & 1+q^3 & q \\ q & q^2 & 1+q^3
\end{array} \right)$}}.
$$
}
\end{example}

\medskip

Cartan matrices provide invariants which are preserved under
derived equivalences and thus improve our understanding
of derived module categories; this is our main motivation
to study normal forms, invariant factors and determinants of Cartan
matrices in this paper.
The following result
is contained in the proof of~\cite[Proposition 1.5]{BS}.

\begin{theorem} \label{unimod}
Let $A$ be a finite-dimensional algebra.
The unimodular equivalence class of the Cartan matrix~$C_A$ is
invariant under derived equivalence.

In particular, the determinant of the Cartan matrix is invariant
under derived equivalence.
\end{theorem}

\begin{remark}
{\rm
We emphasize that the above theorem only deals with
ordinary Cartan matrices $C_A=C_A(1)$.
The determinant
of the $q$-Cartan matrix is in general
not invariant under derived equivalence.
As an example, consider the algebras~$A$ and~$B$ from Example~\ref{brauer-ex},
with $\det C_A(q)=1+q^2+q^4+q^6$ and
$\det C_B(q)=1+q^3+q^6+q^9$.
But, in fact, the algebras~$A$ and~$B$ are derived equivalent;
they are Brauer tree algebras for trees
with the same number of edges and the same exceptional
multiplicity~\cite{R3}.
%In fact, both are Brauer tree algebras (for a definition
%see for instance~\cite{Alperin}) for Brauer trees with
%the same number of edges and the same exceptional multiplicity
%(namely 1 in both cases). Hence, by a result of Rickard
%\cite{R3} they are derived equivalent.
%(Alternatively, one can also show directly
%that $A$ and $B$ are derived equivalent by explicitly constructing
%a suitable tilting complex, as described in the appendix of this paper.)
Note that when specializing $q=1$ we indeed get the same
determinants for the ordinary Cartan matrices, as predicted by
Theorem~\ref{unimod}

However, the natural setting when dealing with $q$-Cartan
matrices is that of graded derived categories.
Indeed, the determinant of the $q$-Cartan matrix (which is  defined
so as to take the grading into account)
is invariant under graded derived equivalences.
We are very grateful to the
referee for pointing
this out to us.
We do not discuss this aspect in this paper further,
but shall address the topic of graded derived equivalences
for gentle algebras in detail in a subsequent publication.

For instance, the above algebra~$B$ is graded derived
equivalent to the algebra~$A$, where the grading on~$A$
is chosen so that~$\alpha$ and~$\beta$ are of degree~2, and
$\delta$ and~$\gamma$ of degree~1. Then Rickard's derived
equivalence~\cite{R3} lifts to a graded derived equivalence.
}
\end{remark}

%%%%%%%%%%%%%%%%%%%%%%%%%%%%%%%%%%%%%%%%%%%%%%%%%%%%%%%%%%%%%%%%%%%%%%%%%%

\section{Gentle algebras}\label{sec:Gentles}

In this section, we shall prove Theorem \ref{Thm1} on
the unimodular equivalence class
of the $q$-Cartan matrix of an arbitrary gentle algebra.
%This extends and refines the main result on the determinant of the
%Cartan matrix of a gentle algebra from~\cite{TH}.

We first recall the definition of special biserial algebras and of gentle
algebras, as these details will be crucial for what follows.
%We then prove a preparatory lemma about $q$-Cartan matrices
%which holds slightly more generally than only for gentle
%algebras.
%The main part of this section will follow, a proof of
%Theorem~\ref{Thm1} on the unimodular equivalence class
%of the $q$-Cartan matrix of an arbitrary gentle algebra.
\smallskip

%Let $K$ be a field and let $Q$ be a quiver, with path algebra
%$KQ$ (see Section~\ref{Sec-def} for background
%on path algebras and quivers with relations).

Let $Q$ be a quiver and $I$ an admissible ideal in the
path algebra~$KQ$.
%(see Section~\ref{Sec-def} for definitions).
%A finite dimensional algebra $A$ is called {\em special
%biserial} if $A$ is Morita equivalent to an algebra $KQ/I$ where
%$I$ is an admissible ideal in~$KQ$ with the following
%properties.
We call the pair $(Q,I)$ a special biserial quiver (with relations)
if it satisfies the following properties.

(i) Each vertex of~$Q$ is starting point of at most two arrows,
and end point of at most two arrows.

(ii) For each arrow~$\alpha$ in~$Q$ there is at most one arrow
$\beta$ such that $\alpha\beta\not\in I$, and at most one arrow
$\gamma$ such that $\gamma\alpha\not\in I$.
\smallskip

A finite-dimensional algebra~$A$ is called special biserial
if it has a presentation as $A=KQ/I$ where $(Q,I)$ is
a special biserial quiver.
\medskip

Gentle quivers form a subclass of the class of
special biserial quivers.
%For each arrow in~$Q$ denote by $s(.)$ its start vertex and by
%$t(.)$ its end vertex.
\smallskip

A pair ($Q,I)$ as above is called a {\em gentle quiver} if it is
special biserial and moreover the following holds.

(iii) The ideal~$I$ is generated by paths of length~2.

(iv) For each arrow~$\alpha$ in~$Q$ there is at most one arrow
$\beta'$ with $t(\alpha)=s(\beta')$ such that $\alpha\beta'\in I$,
and there is at most one arrow~$\gamma'$ with
$t(\gamma')=s(\alpha)$ such that $\gamma'\alpha\in I$.
\smallskip

A finite-dimensional algebra~$A$ is called gentle
if it has a presentation as $A=KQ/I$ where $(Q,I)$ is
a gentle quiver.
\medskip

The following lemma will turn out to be very useful.
It does not only hold for gentle algebras but for those
where we have dropped the final condition (iv) in the definition of
gentle quivers.

%We will study matrices with entries in~$\mathbb{Z}[q]$.
Recall
%the following definition from Section \ref{Sec-intro}.
that
%We call
two matrices $C$, $D$ with entries in~$\mathbb{Z}[q]$
are called unimodularly equivalent (over~$\mathbb{Z}[q]$)
if there exist matrices $P$, $Q$ over~$\mathbb{Z}[q]$ of determinant~1
such that $D=PCQ$.

\begin{lemma}\label{lem:arrow-removal}
Let $(Q,I)$ be a special biserial quiver,
and assume that $I$ is generated by paths of length~2.
Let $A=KQ/I$ be the corresponding special biserial algebra.
Let $\alpha$ be an arrow in~$Q$, not a loop, such that there is no arrow
$\be$ with $s(\al)=t(\be)$ and $\be\alpha\in I$,
or there is no arrow
$\ga$ with $t(\al)=s(\ga)$ and $\al\ga\in I$.
Let $Q'$ be the quiver obtained from~$Q$ by removing
the arrow~$\alpha$, let $I'$
be the corresponding relation ideal
and~$A'=KQ'/I'$.
Then the $q$-Cartan matrices~$C_A(q)$ and~$C_{A'}(q)$ are
unimodularly equivalent
(over~$\Z[q]$).
\end{lemma}

\proof
We consider the case where~$\al$ is an arrow in~$Q$
such that there is no arrow~$\be$
with $s(\al)=t(\be)$ and $\be\al\in I$; the second case is dual.
\\
Let $\al=p_0:v_0\to v_1$. As $(Q,I)$ is special biserial,
there is a unique maximal non-zero path starting with~$p_0$,
say $p=p_0p_1\ldots p_t$, where $p_i:v_i\to v_{i+1}$, $i=1,\ldots,t$.
As $A$ is finite-dimensional, the condition on~$\al=p_0$ guarantees
that $v_i\neq v_0$ for all~$i>0$,
%!!
but we may have $v_i=v_j$
for some $i >j>0$.
Now any non-zero path of length~$j$, say,
ending at~$v_0$ can uniquely be extended to
a non-zero path of length~${j+i}$ ending at~$v_i$,
by concatenation with $p_0\ldots p_{i-1}$.
Conversely, any non-zero path ending at~$v_i$ and involving~$p_0$
arises in this way. \\
Now denote the column corresponding
to a vertex~$v$ in the $q$-Cartan matrix~$C_A(q)$ by~$s_{v}$.
We perform column transformations on~$C_A(q)$
by replacing the columns~$s_{v_i}$ by~$s_{v_i}-q^is_{v_0}$,
for $i=1,\ldots,t+1$
%!!
(if $v_i=v_j$ for some $i > j$, the column $s_{v_i}=s_{v_j}$ will
then be replaced by $s_{v_i}-(q^i+q^j)s_{v_0}$).
The resulting matrix~$\tilde{C}(q)$ is then exactly the Cartan matrix
$C_{A'}(q)$
to the algebra~$A'$ corresponding to the quiver~$Q'$ where $\al=p_0$ has been
removed.
\qed

For any vertex in a quiver~$Q$, its {\em valency} is defined as the number
of arrows attached to it, i.e., the number of incoming arrows plus
the number of outgoing arrows (note that in particular any loop contributes
twice to the valency).
%CB to TH: wieso der referee valency moechte, verstehe ich nicht; degree ist fuer Graphen gaengig

\begin{theorem} \label{Thm:gentle}
Let $(Q,I)$ be a gentle quiver, and let
$A=KQ/I$ be the corresponding gentle algebra.
Denote by $c_k$ the number of minimal
oriented $k$-cycles in~$Q$
with full zero relations.\\
Then the $q$-Cartan matrix~$C_A(q)$
is unimodularly equivalent (over~$\Z[q]$) to a diagonal matrix with entries
$(1-(-q)^k)$, with multiplicity ${c_k}$, $k\geq 1$,
and all further diagonal entries being~$1$.
\end{theorem}

\proof
We want to prove the claim by double induction on the number of vertices
and the number of arrows.
Clearly the result holds if $Q$ has no arrows or if it consists of
one vertex with a loop.

If $Q$ has a vertex~$v$ of valency~1 or~3, or of valency~2 but with
no zero relation at~$v$,
then we can use  Lemma~\ref{lem:arrow-removal}
to remove an arrow from~$Q$;
note that by the conditions in Lemma~\ref{lem:arrow-removal}
the removed arrow is not involved
in any oriented cycle with full zero relations.
Hence $C_A(q)$ is unimodularly equivalent to~$C_{A'}(q)$, where
the corresponding quiver has one arrow less but the same
number of oriented cycles with full zero relations, and
hence the result holds by induction.

Hence we may now assume that all vertices are
of valency 0, 2 or 4, and if a vertex is of valency~2,
then there is a zero relation at the vertex.
Also, if $Q$ is not connected, we may use
induction on the number of vertices to have the result for
the components and thus for the whole quiver; hence
we may assume that $Q$ is connected.
In particular, we now only have vertices $v$ of valency~2 with a non-loop
zero relation at~$v$, and vertices of valency~4.
As we do not have paths of arbitrary lengths, not all
vertices can be of valency~4 (see also~\cite[Lemma~3]{TH}).

Now we take a vertex $v=v_1$ of valency~2, with incoming arrow $p_0:v_0\to v_1$
and outgoing arrow $p_1:v_1\to v_2$ with $p_0p_1=0$
(here, $v_0\neq v\neq v_2$).
\\
As $(Q,I)$ is gentle, there is a unique maximal path~$p$ in~$Q$
with non-repeating arrows
starting in~$v_0$ with~$p_0$,
such that the product of any two consecutive arrows is zero in~$A$;
in our present situation this path is an oriented cycle~$\cal C$ with
full zero relations
returning to~$v_0$.
We denote the vertices on this path by $v_0, v_1=v, v_2, \ldots,
v_s,v_{s+1}=v_0$,
and the arrows by $p_i:v_i\to v_{i+1}$, $i=0,\ldots,s$ (also $p_sp_0=0$);
note that the arrows on~$p$ are distinct, but the vertices are not
necessarily distinct (but we point out that
%!!
$v_i\neq v_1$ for all $i\neq 1$).
\\
Denote by $z_{w}$ the row of the $q$-Cartan matrix~$C_A(q)$
corresponding to the vertex~$w$.
In  $C_A(q)$, we now replace the row~$z_v$ by
the linear combination
$$Z=\sum_{i=1}^{s+1} (-q)^{i-1}z_{v_i} $$
to obtain a new matrix~$\tilde{C}(q)$ (note that this is a unimodular
transformation over~$\Z[q]$).
\\
The careful choice of the coefficients is just made so that we can refine the
argument in~\cite{TH}.  We recall some of the notation there.
For any arrow $\al$ in~$Q$  let $\P(\al)$ be the set of paths
starting with~$\al$
which are non-zero in~$A$.
%!!
At each vertex~$v_i$
there is at most one outgoing arrow $r_i \neq p_i$;
for this arrow we have $p_{i-1}r_i\neq 0$, as $(Q,I)$ is gentle.

Hence, cancelling~$p_i$ induces a natural  bijection
$\phi: \P(p_i) \to \{e_{v_{i+1}}\} \cup \P(r_{i+1})$, for $i=1, \ldots, s-1$,
such that a path of~$q$-weight~$q^j$ is mapped to a path
of~$q$-weight~$q^{j-1}$
%!!
(if there is no arrow~$r_i$, we set $\P(r_i)=\emptyset$).

As $v$ is of valency~2, with a zero-relation at~$v$, we also have
the trivial bijection $\P(p_0)=\{p_0\}\to \{e_{v_{1}}\}$,
again with a weight reduction by~$q$.
Now almost everything cancels in~$Z$,
apart from the one term $1-(-q)^{s+1}$ that we obtain as the entry
in the column corresponding to~$v$.
\\
In the next step, we use the dual (counter-clockwise) operation on the columns
labelled by the vertices on the cycle~${\cal C}$,
i.e., we set $v_{s+2}=v = v_1$ and replace the column~$s_v$  by the
linear combination
$$S=\sum_{i=1}^{s+1} (-q)^{s+1-i}s_{v_{i+1}}\:.$$
Ordering vertices so that
$v$ corresponds to the first row and column of the Cartan matrix,
we have thus unimodularly transformed~$C_A(q)$ to a matrix
of the form
$$\left(
\begin{array}{cccc}
1-(-q)^{s+1} & 0 & \cdots & 0 \\
0 \\
\vdots & & C'(q) \\
0\\
\end{array}
\right)
$$
where $C'(q)$ is the $q$-Cartan matrix of the gentle algebra~$A'$
for the quiver~$Q'$ obtained
from~$Q$ by removing~$v$ and  the arrows incident with~$v$.
Note that in comparison with~$Q$, the quiver~$Q'$ has
one vertex less and one cycle with full zero relations
of length~$s+1$ less; now by induction, the result holds for
$C'(q)=C_{A'}(q)$,
and hence the result for~$C_A(q)$ follows immediately.
\qed

This result has several immediate nice consequences.

\begin{cor} \label{q-det}
Let $(Q,I)$ be a gentle quiver, and let~$A=KQ/I$ be the corresponding
gentle algebra.
Denote by $c_k$ the number of minimal oriented $k$-cycles in~$Q$
with full zero relations.
Then the $q$-Cartan matrix~$C_A(q)$
has determinant $$\det C_A(q)= \prod_{k \geq 1} (1-(-q)^k)^{c_k} \:.$$
\end{cor}

\begin{remark} \label{orient-cycles}
{\rm Let $(Q,I)$ be a gentle quiver, with set of vertices~$Q_0$.
Then, as a direct consequence of Theorem
\ref{Thm:gentle}, there are at most $|Q_0|$ minimal oriented cycles
with full zero relations in the quiver (this could also be proved directly by induction).
}
\end{remark}

Note that the property of being gentle is invariant under
derived equivalence~\cite{SZ},
and we now have some invariants to distinguish the derived equivalence classes.
For a gentle quiver $(Q,I)$,
%and corresponding gentle algebra $A=KQ/I$,
recall that $ec(Q,I)$ and $oc(Q,I)$ denote the number
of minimal oriented cycles in~$Q$ with full zero relations of even and odd
length, respectively.
As an immediate consequence of Corollary~\ref{q-det} we obtain
the following formula for the Cartan determinant which was
the main result in~\cite{TH}:

\begin{cor}\label{det}
Let $(Q,I)$ be a gentle quiver, and let $A=KQ/I$ be the
corresponding gentle algebra.
Then for the determinant of the Cartan
matrix~$C_A$ the following holds.
$$\det C_A = \left\{
\begin{array}{ll} 0 & \mbox{if \: $ec(Q,I)>0$} \\
                  2^{oc(Q,I)} & \mbox{else}
\end{array}
\right.
$$
\end{cor}

Note that in combination with Remark~\ref{orient-cycles}
this implies that the Cartan determinant of a gentle algebra
$A=KQ/I$ is at most~$2^{l(A)}$, where $l(A)=|Q_0|$
is the number of simple modules of~$A$.
\smallskip

The most important application of Theorem~\ref{Thm:gentle}
is the following corollary which gives
for gentle algebras new, combinatorial and easy-to-check invariants
of the derived category.

\begin{cor} \label{ec-oc}
Let $(Q,I)$ and $(Q',I')$ be gentle quivers,
and let $A=KQ/I$ and $A'=KQ'/I'$ be the corresponding
gentle algebras.
If $A$ and~$A'$ are derived equivalent, then
$ec(Q,I)=ec(Q',I')$ and $oc(Q,I)=oc(Q',I')$.
\end{cor}

\proof
Since $A$ and~$A'$ are derived equivalent, their (ordinary)
Cartan matrices $C_A$ and~$C_{A'}$ are unimodularly equivalent
over~$\mathbb{Z}$.
By specializing to $q=1$ in Theorem
\ref{Thm:gentle}, representatives for the
equivalence classes are given by diagonal matrices with
entries~$'2'$ for each minimal oriented cycle with full zero relations of
odd length, an entry~$'0'$ for each such cycle of even length,
and remaining entries~$'1'$.
These are precisely the elementary divisors
over~$\mathbb{Z}$.
% (in fact, they satisfy the necessary divisibility properties).
The elementary divisors of an integer matrix
are uniquely determined, and the diagonal
matrices in Theorem~\ref{Thm:gentle} are actually the
Smith normal forms of $C_A$ and~$C_{A'}$ over~$\mathbb{Z}$.
But by Theorem~\ref{unimod} the unimodular equivalence class,
and hence the Smith normal form,
is invariant under
derived equivalence.

Hence, the diagonal entries in the above normal forms for
$C_A$ and~$C_{A'}$ must occur with exactly the same
multiplicities. Thus we get the same number of minimal oriented
cycles with full zero relations of even length and of odd length,
respectively, i.e., $ec(Q,I)=ec(Q',I')$ and $oc(Q,I)=oc(Q',I')$.
\qed

We now illustrate our results and apply them to derived
equivalence classifications of gentle algebras.

\begin{example} \label{derived-two}
{\bf Gentle algebras with two simple modules.}
{\em There are nine connected gentle quivers $(Q,I)$
with two vertices,
as given in the following list. The dotted lines indicate the
zero relations generating the admissible ideal~$I$.

\vskip1cm

%%%%%%%%%%%%%%%%%%%%%%

\begin{center}
\includegraphics[scale=0.5]{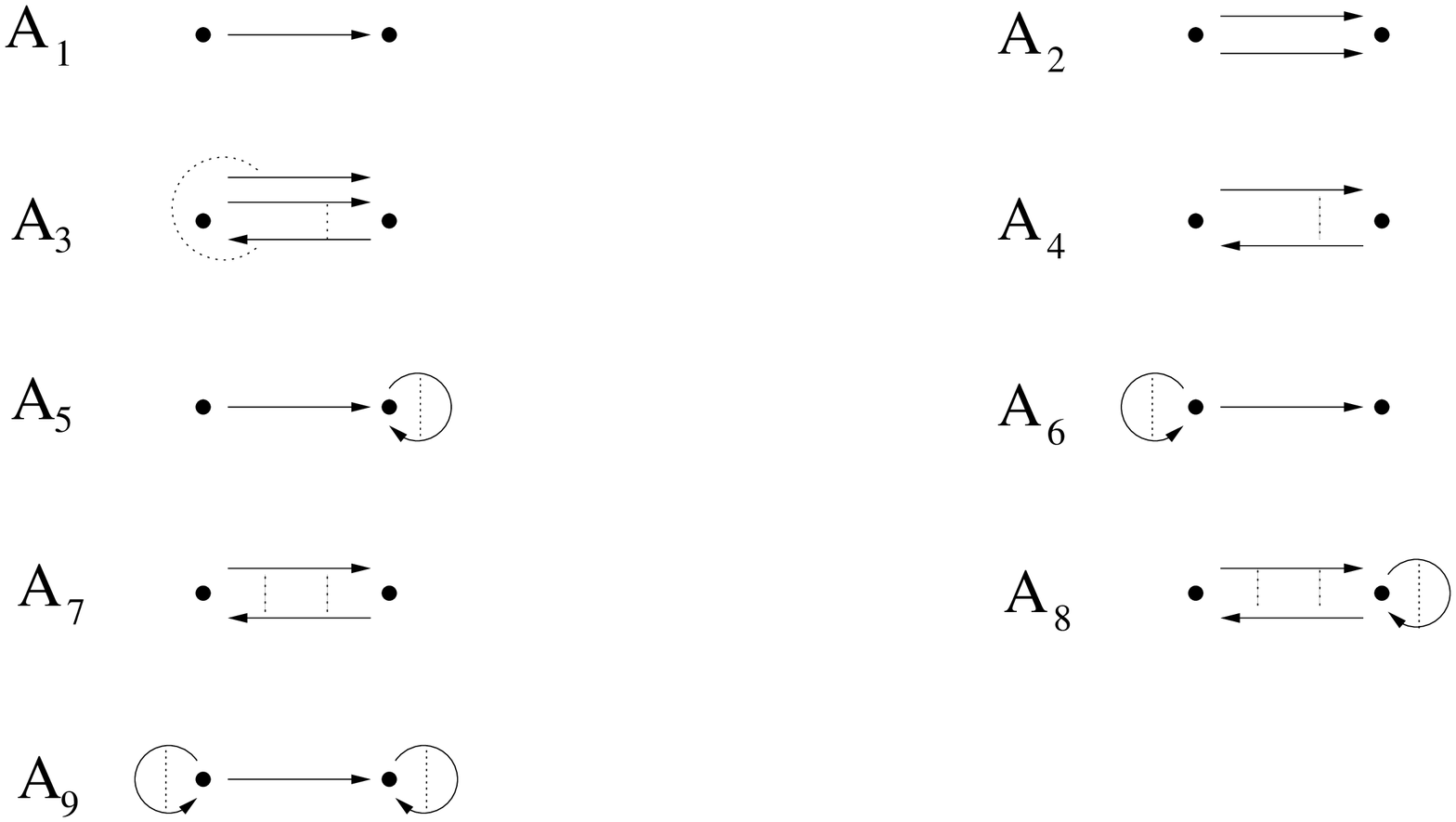}
\end{center}

%%%%%%%%%%%%%%%%%%%%%%

\smallskip

In~\cite{BD} it was shown that these are precisely the basic connected
algebras with two simple modules which are derived tame.
As a direct illustration of our results we show how to classify these
algebras up to derived equivalence.

Recall that the property  of being gentle is invariant under derived
equivalence~\cite{SZ}.
Moreover, the number of simple modules of an algebra is
a derived invariant~\cite{R1}. Thus we will be able to describe
the complete derived equivalence classes.

We have shown above that the numbers $oc(Q,I)$ and $ec(Q,I)$ are derived
invariants. In addition we look at two classical invariants,
the center and the first Hochschild cohomology group~$\HH^1$.
Recall that the center of an algebra (and more generally the
Hochschild cohomology
ring) is invariant under derived equivalence~\cite{R2}.
If the quiver contains a loop, then the dimension of~$\HH^1$
depends on the characteristic being~2 or not. We indicate
the dimension in characteristic~2 in parantheses in the table below.
They can be computed using a method based on work of
M.~Bardzell~\cite{Bardzell}
on minimal projective bimodule resolutions for
monomial algebras; a very nice explicit combinatorial description
is given by C.~Strametz~\cite[Proposition 2.6]{Strametz}.

$$\begin{array}{|c||c||c||c||c||c|c||c||c||c|}
\hline
\mbox{Algebra} & A_1 &  A_2 & A_3 & A_4 & A_5 & A_6 & A_7 & A_8 & A_9 \\
\hline
oc(Q,I) & 0 & 0 & 0 & 0 & 1 & 1 & 0 & 1 & 2 \\
\hline
ec(Q,I) & 0 & 0 & 0 & 0 & 0 & 0 & 1 & 1 & 0 \\
\hline
\dim Z(A) & 1 & 1 & 1 & 2 & 1 & 1 & 1 & 2 & 1 \\
\hline
\dim\HH^1(A) & 0 & 3 & 2 & 1 & 1 (2) & 1 (2) & 1 & 2 (3)& 3 (5) \\
\hline
\end{array}
$$
\smallskip

The algebras $A_1,A_2,A_3,A_4$ are pairwise not derived equivalent.
This can be deduced directly from the above table, since the dimensions
of the first Hochschild cohomology groups are different.

%Alternatively, one can argue as follows
%by using some classical results in tilting theory.
%By~\cite{AH}, the algebras derived equivalent to a hereditary algebra
%of type $\mathbb{A}$ are the gentle algebras whose underlying undirected
%graph is a tree. So $A_1$ can not be derived equivalent to
%any of~$A_2,A_3,A_4$.

%By~\cite{AS}, the algebras derived equivalent to a hereditary algebra
%of type $\tilde{\mathbb{A}}$ are the gentle algebras whose
%underlying graph has exactly one cycle and in this cycle the
%'clock condition' is satisfied (i.e., the number of zero relations
%in clockwise direction and in anti-clockwise direction coincide).
%Hence, $A_2$ is not derived equivalent to $A_3$ or $A_4$.

%Finally, $A_3$ and $A_4$ are not derived equivalent since
%their centers are not isomorphic.
%; in fact, $A_3$ has a one-dimensional
%center and $A_4$ has a two-dimensional center.

The algebras $A_5$ and~$A_6$ are derived equivalent.
(This can be shown by explicitly constructing a suitable tilting
complex, similar to the detailed example given in the Appendix.)
%We leave the details to the reader.)

The algebras $A_7$ and~$A_8$ with Cartan determinant~0 are not
derived equivalent, since their centers have different dimensions.

In summary, there are exactly eight derived equivalence classes of
connected gentle algebras with two simple modules. They are indicated
by the double vertical lines in the above table.
}
\end{example}

\begin{example} \label{derived-three}
{\bf Gentle algebras with three simple modules.}
%{\rm As a second, more involved, illustration
%of our main result we consider (connected)
%gentle algebras with three simple modules.
%Using our results, we
%shall see how easily we can distinguish several gentle algebras
%up to derived equivalence
%which can not be distinguished by the more classical invariants.
%However, to complete the derived equivalence classification we also
%take some classical invariants into account.
{\rm Let $(Q,I)$ be a connected gentle quiver with three vertices,
with corresponding gentle algebra $A=KQ/I$.
By Corollary~\ref{orient-cycles} we deduce that $\det C_A\in \{0,1,2,4,8\}$.
Algebras with different Cartan determinant can not be derived
equivalent, by Theorem~\ref{unimod}.

As an illustration, we shall give a complete derived equivalence
classification of those
algebras with Cartan determinant~0. By Corollary~\ref{det},
a gentle algebra has Cartan determinant~0 if and only
if the quiver contains an even oriented cycle with full zero
relations.
There are 18 connected gentle quivers with three vertices
%simple modules
having Cartan determinant~0, as listed in the following figure.

\vskip1cm

%%%%%%%%%%%%%%%%%%%%%%
\begin{center}
\includegraphics[scale=0.45]{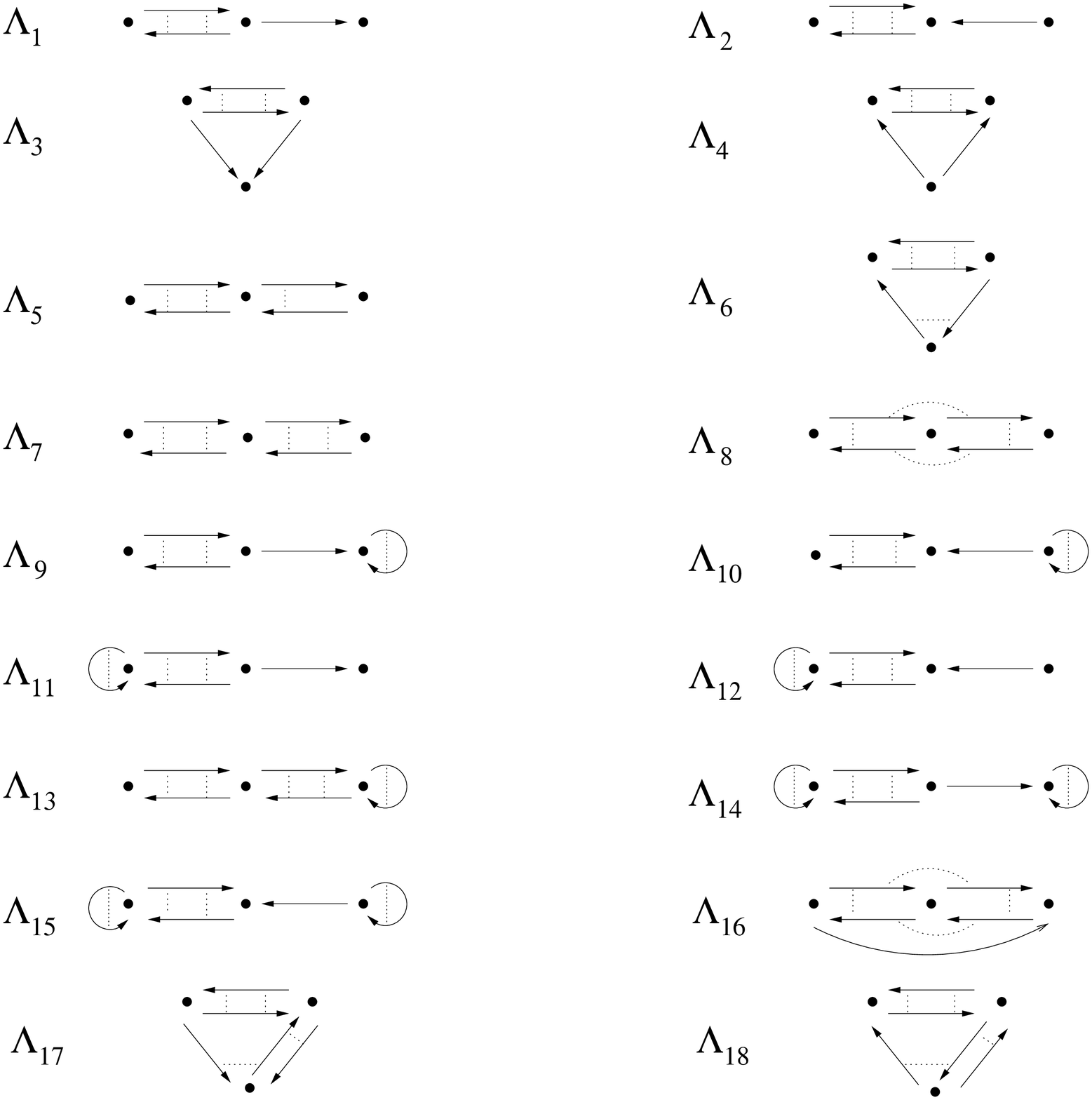}
\end{center}
%%%%%%%%%%%%%%%%%%%%%%

\vskip0.5cm

The main tool will be %our
Corollary~\ref{ec-oc} which states that
the numbers $ec(Q,I)$ and $oc(Q,I)$ of minimal
oriented cycles with full zero relations
of even (resp. odd) length are invariants of the derived category.
This will already settle large parts of the classification.
%especially
%the notoriously difficult problem of distinguishing different
%derived equivalence classes.
In addition we will need to look at the centers and at the first
Hochschild cohomology group. The following
table collects all the necessary invariants.
Again, in the cases where the quiver has loops, the dimension
of~$\HH^1$ depends on the characteristic being~2 or not, and
in these cases the dimension in characteristic~2 is given in parantheses.
\smallskip

$$\begin{array}{|c||c|c|c|c|}
\hline
\mbox{Algebra} & oc(Q,I) & ec(Q,I) & \dim Z(\Lambda) &
\dim\HH^1(\Lambda) \\
\hline
\hline
\Lambda_1 & 0 & 1 & 1 & 1 \\
\hline
\Lambda_2 & 0 & 1 & 1 & 1 \\
\hline
\hline
\Lambda_3 & 0 & 1 & 1 & 4 \\
\hline
\Lambda_4 & 0 & 1 & 1 & 4 \\
\hline
\hline
\Lambda_5 & 0 & 1 & 2 & 2 \\
\hline
\Lambda_6 & 0 & 1 & 2 & 2\\
\hline
\hline
\Lambda_7 & 0 & 2 & 1 & 2\\
\hline
\hline
\Lambda_8 & 0 & 1 & 3 & 2\\
\hline
\hline
\Lambda_9 & 1 & 1 & 1 & 2 (3)\\
\hline
\Lambda_{10} & 1 & 1 & 1 & 2 (3) \\
\hline
\Lambda_{11} & 1 & 1 & 1 & 2 (3) \\
\hline
\Lambda_{12} & 1 & 1 & 1 & 2 (3) \\
\hline
\hline
\Lambda_{13} & 1 & 2 & 2 & 3 (4)\\
\hline
\hline
\Lambda_{14} & 2 & 1 & 1 & 4 (6)\\
\hline
\Lambda_{15} & 2 & 1 & 1 & 4 (6)\\
\hline
\hline
\Lambda_{16} & 0 & 1 & 1 & 4 \\
\hline
\hline
\Lambda_{17} & 0 & 1 & 1 & 3 \\
\hline
\Lambda_{18} &  0 & 1 & 1 & 3 \\
\hline
\end{array}
$$
\smallskip

For the derived equivalence classification, it only remains to
consider those algebras having the same invariants.
In the cases where the algebras are in fact derived equivalent,
we leave out the details of the construction of a suitable tilting complex;
in the appendix a detailed example is provided which serves to indicate
the strategy which also works in all other cases.

The algebras $\Lambda_1$ and~$\Lambda_2$ are derived equivalent.
Moreover,
the algebras $\Lambda_3$ and~$\Lambda_4$ are derived equivalent.

Note that $\Lambda_1$
and~$\Lambda_3$ represent different derived equivalence classes
since their first Hochschild cohomology groups have different dimensions.
%(Alternatively, one can argue that
%$\Lambda_1$ has a discrete derived category (in the sense
%of~\cite{V}, to which we also refer for the definition), whereas
%the derived category of~$\Lambda_3$ is not discrete.)

The algebras $\Lambda_5$ and~$\Lambda_6$ are derived equivalent.
(The details for this case are provided in the appendix.)

Similarly, the algebras $\Lambda_9$, $\Lambda_{10}$, $\Lambda_{11}$
and~$\Lambda_{12}$ are derived equivalent,
the algebras $\Lambda_{14}$ and~$\Lambda_{15}$ are derived
equivalent and
moreover, the algebras $\Lambda_{17}$ and~$n\Lambda_{18}$ are derived
equivalent.

The case of~$\Lambda_{16}$ is more subtle. This algebra has exactly
the same invariants
as the algebras $\Lambda_3$ and~$\Lambda_4$.
However, we claim that $\Lambda_{16}$ is not derived equivalent
to~$\Lambda_4$. In fact, the Lie algebra structures on~$\HH^1$
are not isomorphic. Note that with the Gerstenhaber bracket,
the first Hochschild cohomology becomes a Lie algebra. By a result
of B.\ Keller~\cite{Keller}, this Lie algebra structure on~$\HH^1$ is
invariant under derived equivalence. As mentioned before, by work of
M.~Bardzell~\cite{Bardzell} there is an explicit way of computing
$\HH^1$ for a gentle algebra, and  a nice combinatorial version
due to C.\ Strametz~\cite[Proposition 2.6]{Strametz}
(for the additive structure) and~\cite[Theorem 2.7]{Strametz} (for the Lie algebra
structure). With this method one can compute that the four-dimensional
Lie algebras on~$\HH^1(\Lambda_{16})$ and on~$\HH^1(\Lambda_{4})$
are not isomorphic. In fact,
the Lie algebra center of~$\HH^1(\Lambda_{16})$ is two-dimensional, whereas the
Lie algebra center of~$\HH^1(\Lambda_{4})$ has dimension~1.

This completes the derived equivalence classification of connected
gentle algebras with three simple modules and Cartan determinant~0.
The ten derived equivalence classes are indicated in the above table
by the horizontal double lines.
}
\end{example}

%%%%%%%%%%%%%%%%%%%%%%%%%%%%%%%%%%%%%%%%%%%%%%%%%%%%%%%%%%%%%%%%%%%%%%%%

\section{Skewed-gentle algebras} \label{Sec-skewed}

Skewed-gentle algebras were introduced in~\cite{GdlP};
for the notation and definition we follow here mostly~\cite{BMM},
but we try to explain how the construction works rather than
repeating the technical definition from~\cite{BMM}.

We start with a gentle pair $(Q,I)$.
A set $Sp$ of vertices of the quiver~$Q$ is an admissible set of
{\em special} vertices
if the quiver with relations obtained from~$Q$ by adding loops with
square zero at these
vertices is again gentle; we denote this gentle pair by
$(Q^{\rm\small sp},I^{\rm\small sp})$.
The triple $(Q,Sp,I)$ is then called {\em skewed-gentle}.

We want to point out that the admissibility of the set~$Sp$ of special vertices
is both a local as well as a global condition.
Let $v$ be a vertex in the gentle quiver~$(Q,I)$;
then we can only add a loop at~$v$ if $v$ is of valency~1 or~0
or if it is of valency~2 with a zero relation, but not one coming
from a loop.
Hence only vertices of this type are potential special vertices.
But for the choice of an admissible set of special vertices
we also have to take care of the global condition
that after adding all loops, the
pair $(Q^{\rm\small sp},I^{\rm\small sp})$
still does not have paths
of arbitrary lengths.
\medskip

Given a skewed-gentle triple $(Q,Sp,I)$, we now construct a new quiver with
relations $(\hat{Q}, \hat{I})$ by doubling the special vertices,
introducing arrows to and from these vertices corresponding to the previous
such arrows and replacing a previous zero relation at the vertices by a
mesh relation.\\
More precisely, we proceed as follows.
The non-special (or: ordinary) vertices in~$Q$ are also vertices in
the new quiver;
any arrow between non-special vertices as well as corresponding relations
are also kept.
Any special vertex $v\in Sp$ is replaced by two vertices $v^+$ and
$v^-$ in the new quiver.
An arrow $a$ in~$Q$ from a non-special vertex~$w$ to~$v$ (or from~$v$
to~$w$) will be doubled
to arrows $a^{\pm}:w \to v^{\pm}$ (or $a^{\pm}: v^{\pm} \to w$, resp.)
in the new quiver;
an arrow between two special vertices $v,w$ will correspondingly give
four arrows
between the pairs $v^{\pm}$ and~$w^{\pm}$. We say that these new arrows
lie over the arrow~$a$.
Any relation $ab=0$ where $t(a)=s(b)$ is non-special
gives a corresponding zero relation for paths of length~2 with the same
start and end points lying over~$ab$.
If $v$ is a special vertex of valency~2 in~$Q$, then the corresponding
zero relation at~$v$,
say $ab=0$ with $t(a)=v=s(b)$,
is replaced by mesh commutation relations saying that any two paths of
length~2
lying over~$ab$, having the same start and end points but running over
$v^+$ and~$v^-$, respectively,
coincide in the factor algebra to the new quiver with relations
$(\hat{Q}, \hat{I})$.
\\
We will speak of~$(\hat{Q}, \hat{I})$ as a skewed-gentle
quiver {\em covering} the
gentle pair~$(Q,I)$.
\smallskip

A $K$-algebra is then called {\em skewed-gentle}
if it is Morita equivalent to a factor algebra
$K\hat{Q}/\hat{I}$, where $(\hat{Q},\hat{I})$ comes from a
skewed-gentle triple~$(Q,Sp,I)$ as above.
\medskip

{\bf Remark.}
Let $A=KQ/I$ be gentle.
In a gentle quiver, there is at most one non-zero cyclic path
starting and ending at a given vertex;
hence the diagonal entries in the $q$-Cartan matrix
$C_A(q)$ are $1$ or of the form~$1+q^j$, for some $j\in \N$.\\
If a vertex~$v$ in~$Q$ can be chosen as a special vertex
for a covering skewed-gentle quiver~$\hat Q$, then
the corresponding diagonal entry in~$C_A(q)$ is~$1$,
as otherwise we have paths of arbitrary lengths in~$Q^{\rm\small sp}$;
hence in the corresponding $q$-Cartan matrix for the skewed-gentle
algebra~$\hat A$
we have $\begin{pmatrix}1 & 0 \\ 0 & 1\end{pmatrix}$
on the diagonal for the two split vertices $v^{\pm}$ in~$\hat Q$.
\\

\begin{theorem}\label{Thm:skewed-gentle}
Let $(Q,I)$ be a gentle quiver, and $(\hat{Q},\hat{I})$ a covering
skewed-gentle quiver.
Let $\hat A=K\hat Q/\hat I$ be the corresponding skewed-gentle algebra.
Denote by $c_k$ the number of oriented $k$-cycles in~$(Q,I)$ with full
zero relations.\\
Then the $q$-Cartan matrix~$C_{\hat A}(q)$ is unimodularly equivalent
(over~$\mathbb{Z}[q]$)
to a diagonal matrix with entries
$1-(-q)^k$, with multiplicity ${c_k}$, $k\geq 1$,
and all further diagonal entries being~$1$.
\end{theorem}

\proof
Again, we argue by induction on the number of vertices and arrows.
We let $A=KQ/I$ be the gentle algebra
and $C_A(q)$ the $q$-Cartan matrix as before.
\\
If $Q$ has no arrows, then $\hat Q$ is just obtained by doubling the special
vertices, and this still has no arrows, so the result clearly holds.
\\
If $Q$ has an arrow~$\al$ as in Lemma~\ref{lem:arrow-removal}, with a
non-special
$s(\al)$ in the first case, and a non-special~$t(\al)$ in the second
case, respectively,
then we can argue as in the proof of Lemma~\ref{lem:arrow-removal}
to remove~$\al$.
Let us consider again the situation of the first case, so here
$s(\al)=v_0$ is non-special.

Note that a maximal non-zero path~$p$ starting from~$v_0$ with $\al$ or
$\al^{\pm}$ (if~$v_1$ is special) will end on a non-special vertex
(and hence this
maximal path is unique in~$\hat A$); in general, this path will be
longer than the
one taken in~$A$.

In the column transformations,
we only have to be careful at  doubled vertices on the path~$p$;
here we replace both corresponding columns $\hat s_{v_i^{\pm}}$
of~$C_{\hat A}(q)$
by $\hat s_{v_i^{\pm}}-q^i\hat s_{v_0}$.
This leads to the Cartan matrix for the skewed-gentle algebra where $\al$
or~$\al^{\pm}$, respectively, has been removed from~$\hat Q$, which
is a skewed-gentle
cover for  the quiver obtained from~$Q$ by
deleting~$\al$; then the claim follows by induction.

Now assume $Q$ has a source~$v$ which is special (w.l.o.g.\ the first vertex);
the case of a sink is dual.
Then the $q$-Cartan matrix for~$\hat A$ has the form
$$C_{\hat A}(q)= \left(
\begin{array}{ccccc}
1 & 0 & * & \cdots & * \\
0 & 1 & * & \cdots & * \\
0 & 0 \\
\vdots & \vdots & & \hat C'(q) \\
0 & 0 \\
\end{array}
\right)
\sim
\tilde{{\hat C}}(q)=
\left(
\begin{array}{ccccc}
1 & 0 & 0 & \cdots & 0 \\
0 & 1 & 0 & \cdots & 0 \\
0 & 0 \\
\vdots & \vdots & & \hat C'(q) \\
0 & 0 \\
\end{array}
\right)
$$
where $\hat C'(q)$ is the $q$-Cartan matrix of the
skewed-gentle algebra~$\hat A'$
for the quiver~$\hat Q'$ obtained
from~$\hat Q$ by removing $v^+,v^-$ and the arrows incident with~$v^{\pm}$.
Note that~$\hat Q'$
is the skewed-gentle cover for the quiver
$Q'$ which is obtained from~$Q$
by removing~$v$ and the arrow incident with~$v$,
and the choice $Sp'=Sp\setminus\{v\}$ as the set of
special vertices;
in short, we write this as $\hat Q' = \widehat{Q'}$.
Again, using induction the claim follows immediately.

Thus again, we may now assume that $Q$ has only vertices of valency~2
with a (non-loop) zero relation or vertices of valency~4;
note that any  special vertex in~$Q$ has to be of valency~2.
As before, we may also assume that $Q$ (and hence also~$\hat Q$) are
connected.

If there are no non-special vertices,
or if all non-special vertices are of valency~4, then
$Q^{\rm\small sp}$ is not gentle.
Hence $Q$  has a non-special vertex~$v$ of valency~2 with a zero
relation at~$v$.
Let $p_0:v_0\to v$ be the (unique) incoming arrow.
\\
Again we consider the unique maximal path~$p$ in~$Q$
with non-repeating arrows
starting in~$v_0$ with~$p_0$,
such that the product of any two consecutive arrows is zero in~$A$; as before,
we note that in our current situation
$p$ has to be a cycle ${\cal C}=p_0p_1\ldots p_s$,
where $p_i:v_i\to v_{i+1}$, $i=0,\ldots,s$, and $v_{s+1}=v_0$.
As in the previous situation, we note that the arrows are distinct,
but vertices $v_i\neq v_0$ may be repeated.
\\
For a vertex~$w$ in~$Q$ we denote by $z_{w}$ the row of the $q$-Cartan
matrix~$C_A(q)$ corresponding to~$w$.

If $w$ is non-special, we denote by $\hat z_w$
the corresponding row in the $q$-Cartan matrix $\hat C(q)=C_{\hat A}(q)$.
If $w$ is special, then for the two vertices~$w^{\pm}$
we have two corresponding rows $\hat z_{w^{\pm}}$
in the Cartan matrix~$\hat{C}(q)$,
and we then set $\hat z_w=\hat z_{w^+}+\hat z_{w^-}$.\\[1ex]
Before, we have transformed~$C$ by replacing $z_v$ by
$Z=\sum_{i=1}^{s+1} (-q)^{i-1}z_{v_i}$ and obtained  a matrix~$\tilde{C}(q)$.
We now do a parallel transformation on~$\hat{C}(q)$, that is, we replace
$\hat z_v$ by
$$\hat Z=\sum_{i=1}^{s+1} (-q)^{i-1}\hat z_{v_i}\:,$$
and we obtain a matrix~$\tilde{\hat C}(q)$.
We have to compare the differences and check that everything stays
under control for the induction argument. \\[1ex]
If a vertex $v_i$, $1\leq i \leq s$, is special, note that the doubled
contribution in~$\hat z_{v_i}=\hat z_{v_i^+}+\hat z_{v_i^-}$ is needed
on the one hand for the cancellation with the previous row,
and on the other hand to continue around the cycle~$\cal C$.
As $v$ is non-special and of valency~2 with a zero-relation,
we note that as before, in~$\hat Z$ we only have the contribution
$1-(-q)^{s+1}$ at~$v$.

Following this by the parallel operation to the previous column operation
we then  replace the column~$\hat s_v$  by the linear combination
$$\hat S=\sum_{i=1}^{s+1} (-q)^{s+1-i}\hat s_{v_{i+1}}\:,$$
where we use analogous conventions as before.

With $v$ corresponding to the first row and column of the Cartan matrix,
we have thus unimodularly transformed $\hat{C}(q)$ to a matrix
of the form
$$\left(
\begin{array}{cccc}
1-(-q)^{s+1} & 0 & \cdots & 0 \\
0 \\
\vdots & & \hat C'(q) \\
0\\
\end{array}
\right)
$$
where $\hat C'(q)$ is the Cartan matrix of the skewed-gentle algebra~$\hat A'$
for the quiver~$\hat Q'$ obtained
from~$\hat Q$ by removing~$v$ and  the arrows incident with~$v$.
Note that in fact, $\hat Q'=\widehat{Q'}$ in the notation of our
previous proof,
i.e., as explained earlier, $\hat Q'$
is the skewed-gentle cover for the quiver~$Q'$
and the choice $Sp'=Sp\setminus\{v\}$ as the set of
special vertices.

Thus the result follows by induction.
\qed

\begin{remark}{\rm
By comparing Theorem~\ref{Thm:gentle} and Theorem~\ref{Thm:skewed-gentle}
we observe that
the $q$-Cartan matrix~$C_A(q)$ for the gentle algebra $A$ to~$(Q,I)$,
and the $q$-Cartan matrix~$C_{\hat A}(q)$ for a skewed-gentle cover~$\hat A$
are unimodularly equivalent to diagonal matrices which only differ by adding
as many further~$1$'s on the diagonal as there are special vertices
chosen in~$Q$. In particular, with notation as above,
$$\det C_{\hat A}(q) = \det C_A(q)=
\prod_{k \geq 1} (1-(-q)^k)^{c_k} \:.$$
}
\end{remark}

\medskip

This observation has the following immediate consequence when
specializing to~$q=1$.

\begin{cor}
Let $(Q,I)$ be a gentle quiver, and $(\hat{Q},\hat{I})$ a covering
skewed-gentle quiver.
Then the determinant of the ordinary Cartan matrix
of the skewed-gentle algebra $\hat A=K\hat Q/\hat I$
is the same as the one for the  gentle
algebra  $A=KQ/I$, i.e., $\det C_{\hat A}= \det C_{A}$.
\end{cor}

\begin{remark}
{\rm A gentle algebra and a (proper) skewed-gentle algebra may
have the same $q$-invariants but they cannot be
derived equivalent by~\cite{SZ}, Corollary 1.2.}
%? wieso reicht das ?
\end{remark}

%%%%%%%%%%%%%%%%%%%%%%%%%%%%%%%%%%%%%%%%%%%%%%%%%%%%%%%%%%%%%%%%%%%

\section{Appendix: Tilting complexes and derived equivalences,
a detailed example}

This appendix is aimed at providing enough background on tilting
complexes and explicit computations of their endomorphism rings
so that the interested reader can fill in the details in the
derived equivalence classifications of Examples~\ref{derived-two}
and~\ref{derived-three}. We explained there in detail how to distinguish
derived equivalence classes (since this is the main topic of this paper),
but have been fairly short on indicating why certain algebras
in the lists are actually derived equivalent. In this section
we will go through one example in detail; this will indicate
the main strategy which also works in all other cases.

For an algebra~$A$ denote by $D^b(A)$ the bounded derived
category and by $K^b(P_A)$ the homotopy category of bounded
complexes of finitely generated projective $A$-modules.

Two algebras~$A$ and~$B$ are called derived equivalent if
$D^b(A)$ and~$D^b(B)$ are equivalent as triangulated categories.
By J.\ Rickard's theorem~\cite{R1},
this happens if and only there
exists a tilting complex~$T$ for~$A$ such that the endomorphism
ring $\End_{K^b(P_A)}(T)$ in the homotopy category is
isomorphic to~$B$.
A bounded complex~$T$ of projective $A$-modules is called a
tilting complex if the following conditions are satisfied.

(i) $\Hom_{K^b(A)}(T,T[i])=0$ for $i\neq 0$ (where~$[.]$ denotes
the shift operator)

(ii) $\add(T)$, the full subcategory of~$K^b(P_A)$ consisting of
direct summands of direct sums of copies of~$T$, generates
$K^b(P_A)$ as a triangulated category.
\smallskip

In Example~\ref{derived-three} we stated that the algebras
$\Lambda_5$ and~$\Lambda_6$ are derived equivalent. For the convenience
of the reader we recall the definition of these algebras.

\vskip-.2cm

%%%%%%%%%%%%%%%%%%%%%%

\begin{center}
\includegraphics[scale=0.5]{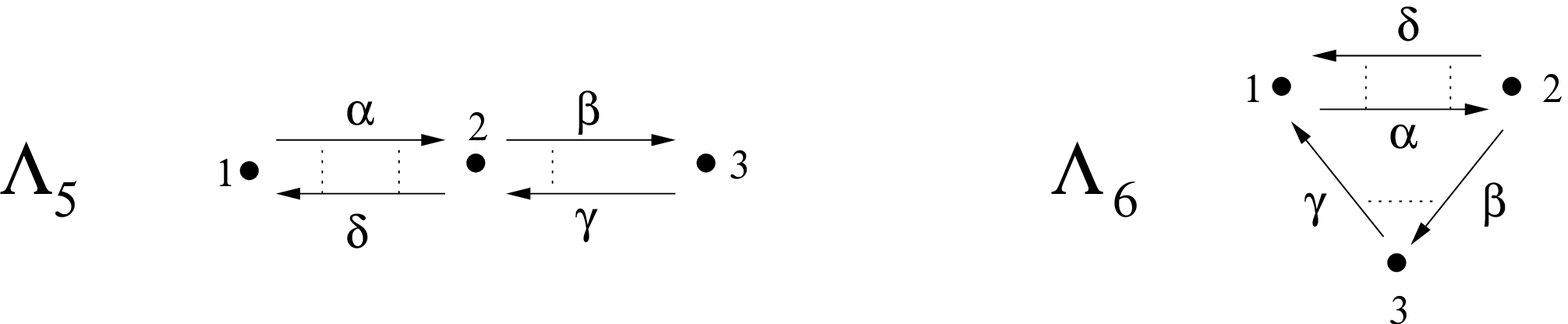}
\end{center}

%%%%%%%%%%%%%%%%%%%%%%

\smallskip

Recall from Section~\ref{Sec-def} our conventions
to deal with right modules and
to read paths from left to right. In particular, left multiplication
by a nonzero path from vertex~$j$ to vertex~$i$ gives a homomorphism
$P_i\to P_j$.
\smallskip

We define the following bounded complex $T:=T_1\oplus T_2\oplus T_3$
of projective $\Lambda_5$-modules. Let
%!!
$T_1\,:\,0\to P_3\to 0$
and $T_3\,:\,0\to P_1\to 0$ be stalk complexes concentrated in degree~0.
Moreover, let $T_2\,:\,0\to P_1\oplus P_3
\stackrel{(\delta,\beta)}{\longrightarrow} P_2\to 0$ (in degrees~0 and~$-1$).
We claim that $T$ is a tilting complex. Property~(i)
above is obvious for all $|i|\ge 2$ since we are dealing with
two-term complexes.

%!!
Let $i=-1$, and consider possible maps $T_2\to T_j[-1]$ where
$j\in\{1,2,3\}$. This is given by a map of complexes as follows
$$\begin{array}{ccccccc}
0 & \longrightarrow & P_1\oplus P_3 &
\stackrel{(\delta,\beta)}{\longrightarrow}
& P_2 & \longrightarrow & 0 \\
 & & & & \downarrow  & & \\
 & & 0 & \longrightarrow & Q & \longrightarrow & \ldots \\
\end{array}
$$
where $Q$ could be either of~$P_1$, $P_3$, or $P_1\oplus P_3$.
But since we are dealing with gentle algebras, no nonzero map
can be zero when composed with both $\delta$ and~$\beta$.
So the only homomorphism of complexes $T_2\to T_j[-1]$ is the
zero map, as desired. Directly from the definition we
see that $\Hom(T_1, T_j[-1])=0$ and $\Hom(T_3, T_j[-1])=0$
(since they are stalk complexes).

Thus we have shown that $\Hom(T, T[-1])=0$.
\smallskip

Now let $i=1$. We have to consider maps $T_j\to T_2[1]$; these are given
as follows
$$\begin{array}{ccccccc}
& & 0 & \longrightarrow & Q & \longrightarrow & \ldots \\
 & & & & \downarrow & & \\
0 & \longrightarrow & P_1\oplus P_3 &
\stackrel{(\delta,\beta)}{\longrightarrow}
& P_2 & \longrightarrow & 0 \\
\end{array}
$$
where $Q$ again can be either of~$P_1$, $P_3$, or $P_1\oplus P_3$.
Now there certainly exist nonzero homomorphisms of complexes.
But they are all homotopic to zero. In fact, every path in the quiver
of~$\Lambda_5$ from vertex~2 to vertex~1 or~3
either starts with~$\delta$ or with~$\beta$. Accordingly, every homomorphism $Q\to P_2$ can be
factored through the map $(\delta,\beta)\,:P_1\oplus P_3\to P_2$.

It follows that $\Hom_{K^b(P_A)}(T, T[1])=0$ (in the homotopy category).
\smallskip

It remains to show that the complex~$T$ also satisfies property (ii)
of the definition of a tilting complex. It suffices to show that
the projective indecomposable modules $P_1$, $P_2$ and~$P_3$,
viewed as stalk complexes, can be generated by~$\add(T)$. This is
clear for $P_1$ and~$P_3$ since they occur as summands of~$T$.
For $P_2$, consider the map of complexes $\Psi:T_2\to T_3\oplus T_1$
given by the identity map on $P_1\oplus P_3$ in degree~0.
Then the stalk complex $P_2[0]$ with~$P_2$ in degree~0 can be shown
to be homotopy equivalent (i.e.\ isomorphic in~$K^b(P_A)$)
to the mapping cone of~$\Psi$. Thus we have a distinguished
triangle
$$\underbrace{T_2}_{\in\add(T)}\to
\underbrace{T_3\oplus T_1}_{\in\add(T)}
\to P_2[0] \to \underbrace{T_2[1]}_{\in\add(T)}.$$
By definition, $\add(T)$ is triangulated, so it follows that also
the stalk complex $P_2[0]\in\add(T)$, which proves (ii).

Hence, $T$ is indeed a tilting complex for~$\Lambda_5$.
\smallskip

By Rickard's theorem, the endomorphism ring of~$T$ in the homotopy
category is derived equivalent to~$\Lambda_5$. We need to show
that $E:=\End_{K^b(P_A)}(T)$ is isomorphic to~$\Lambda_6$.
Note that the vertices of the quiver of~$E$ correspond to the summands
of~$T$.

For explicit calculations, the following formula is useful, which
gives a general
method for computing the Cartan matrix of an endomorphism ring
of a tilting complex from the Cartan matrix of~$A$.
\smallskip

{\em Alternating sum formula.} For a finite-dimensional
algebra~$A$, let $Q=(Q^r)_{r\in\mathbb{Z}}$ and
$R=(R^s)_{s\in\mathbb{Z}}$ be bounded complexes of projective
$A$-modules. Then
$$\sum_i (-1)^i\dim\Hom_{K^b(P_A)}(Q,R[i])=
\sum_{r,s} (-1)^{r-s} \dim\Hom_A(Q^r,R^s).
$$
In particular, if $Q$ and~$R$ are direct summands of a tilting
complex then
$$\dim\Hom_{K^b(P_A)}(Q,R) = \sum_{r,s}(-1)^{r-s}\dim\Hom_A(Q^r,R^s).
$$

Note that the Cartan matrix of~$\Lambda_5$ has the form
\mbox{{\footnotesize
$\left( \begin{array}{ccc} 2 & 2 & 1 \\ 2 & 2 & 1 \\ 1 & 1 & 1
\end{array} \right)$}}. From that,
using the alternating sum formula, we can compute
the Cartan matrix of~$E$ to be
\mbox{{\footnotesize
$\left( \begin{array}{ccc} 1 & 1 & 1 \\ 1 & 1 & 1 \\ 1 & 1 & 2
\end{array} \right)$}}. Note that this is actually the Cartan matrix
of~$\Lambda_6$.

Now we have to define maps of complexes between the summands of~$T$,
corresponding to the arrows of the quiver of~$\Lambda_6$. The final
step then is to show that these maps satisfy the defining relations
of~$\Lambda_6$, up to homotopy.

We define $\tilde{\alpha}:T_1\to T_2$ by the map
$(\alpha\beta,0)\,:\, P_3\to P_1\oplus P_3$
in degree 0. Note that this is indeed a homomorphism of complexes
since $\delta\alpha=0$ in~$\Lambda_5$. Moreover, we define
$\tilde{\beta}:T_2\to T_3$ and $\tilde{\delta}:T_2\to T_1$
by the projection onto the first and second summand in degree 0,
respectively. Finally, we define $\tilde{\gamma}\,:\,T_3\to T_1$
by $\gamma\delta\,:\,P_1\to P_3$.

We now have to check the relations, up to homotopy. We write
compositions from left to right (as in the relations of the
quiver of~$E$).
Clearly,
$\tilde{\alpha}\tilde{\delta}=0$. The composition
$\tilde{\beta}\tilde{\gamma}:T_2\to T_1$ is given in degree 0
by $(\gamma\delta,0)\,:P_1\oplus P_3\to P_3$. So it is not the zero map,
but is homotopic to zero via the homotopy map $\gamma:P_2\to P_3$
(use that $\gamma\beta=0$ in~$\Lambda_5$). Finally, consider
$\tilde{\delta}\tilde{\alpha}$ on~$T_2$. It is given by
$\mbox{{\footnotesize $\left( \begin{array}{cc} 0 & \alpha\beta \\ 0 & 0
\end{array} \right)$}}$ in degree 0 and the zero map in degree~$-1$.
It is indeed homotopic to zero via the homotopy map
$(\alpha,0)\,:\,P_2\to P_1\oplus P_3$. (Note that
here we use that $\alpha\delta=0$ and $\delta\alpha=0$
in~$\Lambda_5$.)

Thus, we have defined maps between the summands of~$T$, corresponding
to the arrows of the quiver of~$\Lambda_6$. We have shown that they
satisfy the defining relations of~$\Lambda_6$, and that the Cartan
matrices of~$E$ and~$\Lambda_6$
%!!
coincide.
From this we can conclude that
$E\cong \Lambda_6$. Hence, $\Lambda_5$ and
$\Lambda_6$ are derived equivalent, as desired.
\smallskip

All the other derived equivalences stated in Examples~\ref{derived-two}
and~\ref{derived-three} can be verified exactly along these lines.
In particular, they can also be realized by tilting complexes
with non-zero entries in only two degrees.

%%%%%%%%%%%%%%%%%%%%%%%%%%%%%%%%%%%%%%%%%%%%%%%%%%%%%%%%%%%%%%%%%%%%%


\begin{thebibliography}{99}

%\bibitem{Alperin}
%J. Alperin, Local representation theory,
%Cambridge Studies in Advanced Mathematics, 11. Cambridge University
%Press, Cambridge, 1986

\bibitem{AH}
I. Assem, D. Happel,
Generalized tilted algebras of type~$\mathbb{A}_n$,
Comm. Algebra 9 (1981), no.20, 2101-2125

\bibitem{AS}
I.\ Assem, A.\ Skowro\'{n}ski,
Iterated tilted algebras of type~$\tilde{\mathbb{A}}_n$,
Math.\ Z.\ 195 (1987), 269-290

\bibitem{Bardzell}
M.~J.\ Bardzell,
The alternating syzygy behavior of monomial algebras,
J.\ Algebra 188 (1997) 69--89

\bibitem{BD}
V.\ Bekkert, Y.\ Drozd,
Derived categories for some classes of algebras,
Preprint. Available at
{\tt http://www.mat.ufmg.br/$\sim$bekkert/pub.html}

\bibitem{BMM}
V. Bekkert, E.N. Marcos, H.A. Merklen,
Indecomposables in derived categories of skewed-gentle algebras,
Communications in Algebra (6) 31 (2003) 2615-2654

\bibitem{BGS}
G. Bobi\'{n}ski, C. Geiss, A. Skowro\'{n}ski,
Classification of discrete derived categories,
Cent. Eur. J. Math. 2 (2004)
%, no. 1,
19-49

\bibitem{BS}
R.\ Bocian, A.\ Skowro\'{n}ski,
Weakly symmetric algebras of Euclidean type,
J. Reine Angew. Math. 580 (2005) 157-200


%\bibitem{Fa} R. Farnsteiner,
%Trivial extensions, Notes, Jan.~7, 2005

\bibitem{F}
K.\ Fuller, The Cartan determinant and global dimension of artinian rings,
in: Azumaya algebras, actions, and modules, Proceedings 1990,
Contemporary Mathematics 124 (AMS 1992) 51-72

\bibitem{GK}
C.\ Geiss, H.\ Krause,
On the notion of derived tameness,
J. Algebra Appl. 1 (2002) %no. 2
133-157

\bibitem{GdlP}
C.\ Geiss, J.~A.\ de la Pe$\tilde{\textrm n}$a,
Auslander-Reiten components for clans,
Boll. Soc. Mat. Mexicana 5 (1999) 307-326

\bibitem{TH}
T.\ Holm, Cartan determinants for gentle algebras,
Archiv Math. 85 (2005), 233-239

\bibitem{Keller}
B.\ Keller,
Derived invariance of higher structures on the Hochschild complex,
Preprint (2003)

\bibitem{R1}
J.\ Rickard,
Morita theory for derived categories,
J.\ London Math.\ Soc.\ (2) 39  (1989) 436-456

\bibitem{R2}
J.\ Rickard,
Derived equivalences as derived functors.
J.\ London Math.\ Soc.\ (2) 43 (1991) 37-48

\bibitem{R3}
J.\ Rickard,
Derived categories and stable equivalence.
J.\ Pure Appl.\ Algebra 61 (1989) 303-317

\bibitem{SZ}
J.\ Schr\"oer, A.\ Zimmermann,
Stable endomorphism algebras of modules over special biserial algebras,
Math.\ Z.\ 244 (2003) 515-530

\bibitem{Strametz}
C.\ Strametz, The Lie algebra structure of the first Hochschild
cohomology group of a monomial algebra,
Preprint (2001), math.RT/0111060; announced in:
C.\ R.\ Math.\ Acad,\ Sci.\ Paris 334 (2002) 733-738

\bibitem{V}
D.\ Vossieck,
The algebras with discrete derived category,
J.\ Algebra 243 (2001) 168-176


\end{thebibliography}
\end{document}